\documentclass[12pt]{amsart}

\input xy
\usepackage[english]{babel}
\usepackage[latin1]{inputenc}
\usepackage{graphicx}
\usepackage{amsmath}
\usepackage{amssymb}
\usepackage{amscd}
\usepackage{color}
\usepackage{oldgerm}
\usepackage{amsfonts}
\usepackage{newlfont}
\usepackage{longtable}
\usepackage{multirow}
\usepackage{amsthm}

\theoremstyle{plain}\swapnumbers
\usepackage[all]{xy}

\usepackage{upref}

\def\Ker{\operatorname{Ker}}

\def\lim{\operatorname{lim}}

\DeclareOldFontCommand{\rm}{\normalfont\rmfamily}{\mathrm}

\theoremstyle{plain}\swapnumbers

\newtheorem{Theorem}{Theorem}[section]

\newtheorem{Prop}[Theorem]{Proposition}

\newtheorem{claim}{Claim}
\newtheorem{prop A}{Proposition A}

\title{Covariant functors commuting with direct limits}

\author{R. Garc\'ia-Delgado}

\address{Centro de Investigaci\'on en Matem\'aticas A. C - Unidad M\'erida, Yucat\'an, M\'exico; Carretera Sierra Papacal Chuburna Puerto Km 5, 97302, Yucat\'an, M\'exico.}

\email{rosendo.garciadelgado@alumnos.uaslp.edu.mx}



\date{\today}

\begin{document}

\maketitle

\begin{abstract}
In this work we state conditions on a covariant right exact functor so that it commutes with direct limits. These conditions are related to the commutativity of the functor under direct limits of projective modules. We prove that if the functor commutes with direct limits of projective modules, then the functor commutes with direct limits. 
\end{abstract}

\section*{Introduction}

The main examples of functors that commute with direct and inverse limits are $\otimes_{\Lambda}$ and $\operatorname{Hom}_{\Lambda}$, as well as the left derived functors $\operatorname{Tor}^{\lambda}_n$ (see \cite{Cartan}, {\bf Prop 9.2}, chapter V, \S 9 and {\bf Prop. 1.3}, chapter VI, \S 1). The commutative properties in functors give information about the structure module.  In \cite{Lenzing} is proved that a right $\Lambda$-module $M$ is finitely presented if and only if the functor $\operatorname{Hom}_{\Lambda}(M,\cdot)$, preserves direct limits. Similarly, in \cite{Breaz} is proved that $\operatorname{Ext}_{\Lambda}^1(M,\cdot)$ commutes with inverse limits if and only if $M$ is projective.
\smallskip

In {\bf Thm. \ref{teorema 2}}, we state conditions in a covariant right exact functor $T$, so that it commutes with direct limits based on the commutativity of $T$ under direct limits of projective modules. 

\section{Direct limits of modules}

We start by reviewing some definitions about direct system of modules, as well as some basic results and the notation that will be used in the sequel. Throughout this work we consider \emph{left} $\Lambda$-modules.
\smallskip

Let $\Lambda$ be a ring and let $\{A^{\gamma}\mid \gamma \in \Omega\}$ be an indexed family of $\Lambda$-modules by a set $\Omega$. The {\bf direct product} {\boldmath $\prod_{\gamma \in \Omega}A^{\gamma}$}, is the set consisting of those maps $f:\Omega \to \cup_{\gamma}A^{\gamma}$ such that $f(\gamma)$ belongs to $A^{\gamma}$ for each $\gamma$ in $\Omega$. The direct product $\prod_{\gamma \in \Omega}A^{\gamma}$ is a left $\Lambda$-module with coordinatewise addition and usual left scalar multiplication. 
\smallskip

The {\bf direct sum} {\boldmath $\oplus_{\gamma}A^{\gamma}$} is the $\Lambda$-submodule of $\prod_{\gamma \in \Omega}A^{\gamma}$, consisting of those maps $f$ in $\prod_{\gamma \in \Omega}A^{\gamma}$ such that $f(\gamma)=0$ except in a finite number of $\gamma$'s. 
\smallskip

Any element $f$ in $\prod_{\gamma \in \Omega}A^{\gamma}$ can be written as $f=\left(x_{\gamma}\right)_{\gamma \in \Omega}$, where $f(\gamma)=x_{\gamma}$ is in $A^{\gamma}$, and $\gamma$ is in $\Omega$. For each $\alpha$ in $\Omega$ there is a map $p_{\alpha}:\prod_{\gamma \in \Omega}A^{\gamma} \to A^{\alpha}$, given by $p_{\alpha}(f)=x_{\alpha}$. Using the Kronecker delta map $\delta_{\alpha,\beta}$, each $\alpha$ in $\Omega$ defines a map $c^{\alpha}:A^{\alpha} \to \oplus_{\gamma}A^{\gamma}$ as follows: Let $x$ be in $A^{\alpha}$, then $c^{\alpha}(x):\Omega \to \cup_{\gamma}A^{\gamma}$ is defined by $c^{\alpha}(x)(\beta)=\delta_{\alpha,\beta}\,x$. The map $c^{\alpha}$ is known as the \textbf{$\alpha$-th injection}. 
\smallskip

Then any $f=(x_{\gamma})_{\gamma \in \Omega }$ in $\bigoplus_{\gamma \in \Omega}A^{\gamma}$, can be written as a finite sum of the form $f=\displaystyle{\sum_{\gamma \in \Omega}}c^{\gamma}\left(x_{\gamma}\right)$.
\smallskip

Let $\Omega$ be an ordered partially set. A {\bf direct system} of $\Lambda$-modules over $\Omega$ is an ordered pair $\{A^{\alpha}; \varphi_{\alpha}^{\beta}\}$, consisting of an index family of $\Lambda$-modules $\{A^{\gamma}\mid \gamma \in \Omega\}$, together with a family of homomorphisms $\{\varphi_{\alpha}^{\beta}:A^{\alpha} \to A^{\beta}\mid \alpha \leq \beta\}$, such that $\varphi_{\alpha}^{\alpha}=\operatorname{Id}_{A_{\alpha}}$ and $\varphi^{\gamma}_{\alpha}=\varphi^{\gamma}_{\beta} \circ \varphi^{\beta}_{\alpha}$, for all $\alpha \leq \beta \leq \gamma$.
\smallskip

The {\bf direct limit} of $\{A^{\alpha}; \varphi^{\alpha}_{\beta}\}$ is a $\Lambda$-module $\underrightarrow{\operatorname{lim}}A^{\gamma}$ and a family of $\Lambda$-homomorphisms $\{\sigma^{\alpha}:A^{\alpha} \rightarrow \underrightarrow{\operatorname{lim}}A^{\gamma}\,\mid \, \alpha \in \Omega \}$, such that:

\textbf{(i)} $\sigma^{\beta} \circ \varphi^{\alpha}_{\beta}=\sigma^{\alpha}$, whenever $\alpha \leq \beta$.
\smallskip

\textbf{(ii)} For every $\Lambda$-module $X$ with $\Lambda$-homomorphisms $\psi^{\gamma}:A^{\gamma} \to X$ satisfying $\psi^{\beta} \circ \varphi^{\alpha}_{\beta}=\psi^{\alpha}$, for all $\alpha \leq \beta$, there exists a unique $\Lambda$-homomorphism $\theta\!:\! \!\underrightarrow{\operatorname{lim}}A^{\gamma} \to X$, such that $\theta \!\circ\! \sigma^{\gamma}\!\!=\!\!\psi^{\gamma}$, for all $\gamma$ in $\Omega$. We shall prove that the direct limits of modules exists. Let $\bar{\theta}$ be the map between $\oplus_{\gamma \in \Omega}A^{\gamma}$ and $X$ defined by:
$$
\bar{\theta}(f)=\sum_{\gamma \in \Omega}\psi^{\gamma}(x_{\gamma}),\,\,\,\text{ where }f=\sum_{\gamma \in \Omega}c^{\gamma}(x_{\gamma}).
$$
Let $N$ be the $\Lambda$-submodule of $\oplus_{\gamma \in \Omega}A^{\gamma}$, given by:
\begin{equation}\label{submodulo N}
N=\operatorname{Span}_{\Lambda}\{c^{\beta}\left(\varphi^{\alpha}_{\beta}(x)\right)-c^{\alpha}(x)\mid x \in A^{\alpha} \text{ and }\alpha \leq \beta \}.
\end{equation}
Then $N$ is contained in the kernel of $\bar{\theta}$. Thus, there is a $\Lambda$-homomorphism $\theta$ between the quotient $\oplus_{\gamma \in \Omega}A^{\gamma}/N$ and $X$ such that $\theta \circ \sigma^{\alpha}=\psi^{\alpha}$, for all $\alpha$ in $\Omega$, where $\sigma^{\alpha}$ is the map defined by:
\begin{equation}\label{canonical map of direct limits}
\begin{array}{rccl}
\sigma^{\alpha}:& A^{\alpha} & \longrightarrow & \oplus_{\gamma \in \Omega}A^{\gamma}/N\\
& x & \mapsto & c^{\alpha}(x)+N,\,\,\text{ for all }\,\alpha \in \Omega,
\end{array}
\end{equation}
Then $\underrightarrow{\operatorname{lim}}A^{\gamma}$ exists and is equal to $\oplus_{\gamma \in \Omega}A^{\gamma}/N$. The map $\sigma^{\alpha}$, is the {\bf canonical map of} {\boldmath $\underrightarrow{\operatorname{lim}}A^{\gamma}$} and sometimes we write it as $\sigma^{\alpha}_A$.
\smallskip

A partially ordered set $\Omega$ is a {\bf directed set} if for any $\alpha,\beta$ in $\Omega$, there exists $\gamma$ in $\Omega$ such that $\gamma \geq \alpha$ and $\gamma \geq \beta$.
There are two reasons to consider direct systems over directed sets. The first one is that a simpler description of the elements in the direct limit can be given; the second one is that $\underrightarrow{\operatorname{lim}}$ preserves short exact sequences. 

\begin{Prop}[\cite{Rotman}, \textbf{Prop. 7.98}]\label{rotman}{\sl
Let $\{A^{\alpha};\varphi^{\alpha}_{\beta}\}$ be a direct system of $\Lambda$-modules over a directed set $\Omega$, and let $c^{\alpha}:A^{\alpha} \to \oplus_{\gamma}A^{\gamma}$ be the $\alpha$-th injection and let $\underrightarrow{\operatorname{lim}}A^{\gamma}=\oplus_{\gamma \in \Omega}A^{\gamma}/N$, where $N$ is the submdole defined in \eqref{submodulo N}. Then,

\textbf{(i)} each element of $\underrightarrow{\operatorname{lim}}A^{\gamma}$, has a representative of the form $c^{\gamma}(a)+N$, for some $a$ in $A^{\gamma}$ and $\gamma$ in $\Omega$.
\smallskip

\textbf{(ii)} $c^{\gamma}(a)+N=0$ if and only if $\varphi^{\alpha}_{\beta}(a)=0$, for some $\beta \geq \alpha$.
}
\end{Prop}
From now on we will consider direct systems over a unique directed index set $\Omega$. In addition $\alpha,\beta$ and $\gamma$ always will denote elements in $\Omega$. 
\smallskip

A {\bf homomorphism} between the direct systems $\{C^{\alpha};\phi^{\alpha}_{\beta}\}$ and $\{A^{\alpha};\varphi^{\alpha}_{\beta}\}$, is a family of $\Lambda$-homomorphisms $\{f^{\gamma}:C^{\gamma} \to A^{\gamma} \mid \gamma \in \Omega \}$, such that $f^{\beta} \circ \phi^{\alpha}_{\beta}=\varphi^{\alpha}_{\beta} \circ f^{\alpha}$, for all $\alpha \leq \beta$. If there exists such a homomorphism, then there exists a $\Lambda$-homomorphism $F:\underrightarrow{\operatorname{lim}}C^{\gamma} \to \underrightarrow{\operatorname{lim}}A^{\gamma}$, verifying $F \circ {\sigma_C}^{\alpha}=\sigma_A^{\alpha} \circ f^{\alpha}$, for all $\alpha$, where ${\sigma_C}^{\alpha}:C^{\alpha} \to \underrightarrow{\operatorname{lim}}C^{\gamma}$ (resp. $\sigma_A^{\gamma}:A^{\alpha} \to \underrightarrow{\operatorname{lim}}A^{\gamma}$) is the canonical map of $\underrightarrow{\operatorname{lim}}C^{\gamma}$ (resp. $\underrightarrow{\operatorname{lim}}A^{\gamma}$).
\smallskip

Let $\{C^{\alpha};\phi^{\alpha}_{\beta}\}$, $\{A^{\alpha};\varphi^{\alpha}_{\beta} \}$ and $\{E^{\alpha};\kappa^{\alpha}_{\beta}\}$, be direct systems such that:
$$
\xymatrix{
0 \ar[r] & C^{\alpha} \ar[r]^{f^{\alpha}} & A^{\alpha} \ar[r]^{g^{\alpha}} & E^{\alpha} \ar[r] & 0
},\quad \text{ where } \alpha \in \Omega,
$$
is a short exact sequence in the category of direct system of modules. Then the following diagram is commutative with rows exacts:
\begin{equation}\label{prop exact lim directo}
\xymatrix{
0 \ar[r] & C^{\alpha} \ar[r]^{f^{\alpha}} \ar[d]_{\sigma_C^{\alpha}} & A^{\alpha} \ar[r]^{g^{\alpha}} \ar[d]_{\sigma_A^{\alpha}} & E^{\alpha} \ar[r] \ar[d]_{\sigma_E^{\alpha}} & 0\\
\, & \underrightarrow{\operatorname{lim}}C^{\gamma} \ar[r]^{F} & \underrightarrow{\operatorname{lim}}A^{\gamma} \ar[r]^{G} & \underrightarrow{\operatorname{lim}}E^{\gamma} \ar[r] & 0
}
\end{equation}
Let $T:\operatorname{Mod}_{\Lambda} \to \operatorname{Mod}_{\Lambda_1}$ be a covariant functor. Then $\{T\left(A^{\alpha}\right);T\left(\varphi_{\beta}^{\alpha}\right)\}$ is a direct system of $\Lambda_1$-modules. Observe that there exists a $\Lambda_1$-homomorphism induced by the canonical map $\sigma_A^{\gamma}$,
\begin{equation}\label{Lsigma}
\begin{array}{rccl}
\sigma^T_A:& \underrightarrow{\operatorname{lim}}T\left( A^{\gamma} \right) & \longrightarrow & T\left(\underrightarrow{\operatorname{lim}}A^{\gamma}\right)\\
\,& c^{\gamma}(x)+N & \mapsto & T\left(\sigma_A^{\gamma}\right)(x),\,\,\text{ where }\,\gamma \in \Omega\,\text{ and }x \in T(A^{\gamma}).
\end{array}
\end{equation}
The functor $T$ is {\bf of the type} {\boldmath $L\Sigma^{*}$} if $\sigma^T_A$ is an \emph{isomorphism}; that is, if $T$ \emph{commutes with direct limits}.

\section{Functors commuting with direct limits}

Any $\Lambda$-module $A$ can be embedded into an exact sequence of the form:
$$
\xymatrix{ 
 0 \ar[r] &
M \ar[r] & 
P \ar[r] &
A\ar[r] &
 0 }
$$
with $P$ a projective module (see \cite{Cartan}, {\bf Thm. 2.3}, Chapter I \S2). Now we are in conditions to state the main result of this work.

\begin{Theorem}\label{teorema 2}{\sl
Let $T$ be a one variable covariant half exact functor. If $T$ commutes with direct limits of projective modules, then $T$ is of the type $L\Sigma^{*}$.}
\end{Theorem}
\begin{proof}
Let $\{A^{\alpha};\varphi^{\alpha}_{\beta}\}$ be a direct system of $\Lambda$-modules and let $A=\underrightarrow{\operatorname{lim}}A^{\gamma}$. We shall prove that $\sigma^T_A$ (see \eqref{Lsigma}) is an isomorphism. First we will prove the following result.

\begin{claim}{\sl
For each $\gamma$ in $\Omega$, there exists a short exact sequence:
$$
\xymatrix{
0 \ar[r] & M^{\gamma} \ar[r] & P^{\gamma} \ar[r] & A^{\gamma} \ar[r] & 0
}
$$
where $P^{\gamma}$ is a free $\Lambda$-module, and $\{P^{\alpha};\phi^{\alpha}_{\beta}\}$ and $\{M^{\alpha};\varphi^{\alpha}_{\beta}\}$ are direct system of $\Lambda$-modules.
}
\end{claim}
Indeed, let $P^{\gamma}$ be the set consisting of those maps $f:A^{\gamma} \to \Lambda$ such that $f(a)=0$, except for a finite number of $a$'s. Then $P$ is a $\Lambda$-module with the usual operations.   
\smallskip

For each $a$ in $A^{\gamma}$, let $e^{\gamma}(a)$ in $P^{\gamma}$ be defined by $e^{\gamma}(a)(b)=\delta_{a,b}$, for all $b$ in $A^{\gamma}$. Then any $f$ in $P^{\gamma}$ can be written as: $f=\displaystyle{\sum_{a \in A^{\gamma}}}f(a)e^{\gamma}(a)$. This shows that $P^{\gamma}$ is the free $\Lambda$-module generated by $e^{\gamma}(a)$, where $a$ is in $A^{\gamma}$. We extend the map $\varphi_{\beta}^{\alpha}:A^{\alpha} \to A^{\beta}$, to $P^{\alpha}$ and $P^{\beta}$, $\alpha \leq \beta$, by:
$$
\begin{array}{rccl}
\phi_{\beta}^{\alpha}:& P^{\alpha} & \longrightarrow & P^{\beta}\\
& f & \mapsto &  \displaystyle{\sum_{a \in A^{\alpha}}}f(a)e^{\beta}\left(\varphi^{\alpha}_{\beta}(a)\right).
\end{array}
$$
Then $\phi_{\gamma}^{\beta} \circ \phi_{\beta}^{\alpha}=\phi_{\gamma}^{\alpha}$, for all $\alpha \leq \beta \leq \gamma$, and $\{P^{\alpha}; \phi_{\alpha}^{\beta} \}$ is a direct system of free $\Lambda$-modules. Let $q^{\gamma}:P^{\gamma} \to A^{\gamma}$ be the $\Lambda$-homomorphism defined by $q^{\gamma}\left(e^{\gamma}(a)\right)=a$, for all $a$ in $A^{\gamma}$. Let $M^{\gamma}=\Ker(q^{\gamma})$ and let $\iota^{\gamma}:M^{\gamma} \to P^{\gamma}$ be the inclusion map. We have the short exact sequence of $\Lambda$-modules:
$$
\xymatrix{ 
 0 \ar[r] &
M^{\gamma} \ar[r]^{\iota^{\gamma}} & 
P^{\gamma} \ar[r]^{q^{\gamma}} &
A^{\gamma}\ar[r] &
 0 }
$$
Also we have $q^{\beta}\circ \phi^{\alpha}_{\beta}=\varphi^{\alpha}_{\beta}\circ q^{\alpha}$, where $\alpha \leq \beta$., which implies that $\phi_{\beta}^{\alpha}\left(M^{\alpha}\right) \subset M^{\beta}$. Let $\chi_{\beta}^{\alpha}=\phi_{\beta}^{\alpha}|_{M^{\alpha}}:M^{\alpha} \to M^{\beta}$, for $\alpha \leq \beta$. Then $\{M^{\alpha};\chi_{\beta}^{\alpha}\}$ is a direct system of $\Lambda$-modules. In addition we obtain the commutative diagram with exact rows:
$$
\xymatrix{ 
0 \ar[r] & M^{\alpha} \ar[d]_{\chi_{\beta}^{\alpha}} \ar[r]^{\iota^{\alpha}} & P^{\alpha} \ar[d]_{\phi_{\beta}^{\alpha}} \ar[r]^{q^{\alpha}} & A^{\alpha} \ar[d]_{\varphi_{\beta}^{\alpha}} \ar[r] & 0\\
0 \ar[r] & M^{\beta}  \ar[r]^{\iota^{\beta}} & P^{\beta} \ar[r]^{q^{\beta}} & A^{\beta} \ar[r] & 0
}
$$
\begin{claim}{\sl
The short sequence:
$$
\xymatrix{
0 \ar[r] & \underrightarrow{\operatorname{lim}}M^{\gamma} \ar[r] & \underrightarrow{\operatorname{lim}}P^{\gamma} \ar[r] & \underrightarrow{\operatorname{lim}}A^{\gamma}\ar[r] & 0
}
$$
is exact and $\underrightarrow{\operatorname{lim}}P^{\gamma}$ is a free $\Lambda$-module.}
\end{claim}
Indeed, let $P$ be $\Lambda$-module consisting of maps $f:A \to \Lambda$ such that $f(x)=0$ except for a finite number of $x$'s. For each $x$ in $A$, let $e(x)$ be the element of $P$ defined by $e(x)(y)=\delta_{x,y}$, for all $y$ in $A$. Then any $f$ in $P$ can be written as $f=\displaystyle{\sum_{x \in A}}f(x)e(x)$, and $P$ is a free $\Lambda$-module generated by $e(x)$, where $x$ is in $A$. 
\smallskip

Let $q:P \to A$ be the $\Lambda$-homomorphism defined by $q(e(x))=x$, for all $x$ in $A$. Let $M=\Ker(q) \subset P$ and let $\iota:M \to P$ be the inclusion map. Then we have the following short exact sequence of $\Lambda$-modules:
$$
\xymatrix{ 
 0 \ar[r] &
M \ar[r]^{\iota} & 
P \ar[r]^{q} &
A\ar[r] &
 0 }
$$
For each $\gamma$, let $\tau^{\gamma}:P^{\gamma} \to P$ be defined by $\tau^{\gamma}\left(e^{\gamma}(a)\right)=e\left(\sigma_A^{\gamma}(a)\right)$, for all $a$ in $A^{\gamma}$. 
Then $\tau^{\beta} \circ \phi_{\beta}^{\alpha}=\tau^{\alpha}$, for all $\alpha \leq \beta$. Let $\sigma^{\alpha}_P:P^{\alpha} \to \underrightarrow{\operatorname{lim}}P^{\gamma}$ be the canonical map. Each of the $\Lambda$-homomorphisms $\tau^{\gamma}$, yields a $\Lambda$-homomorphism $\tau:\underrightarrow{\operatorname{lim}}P^{\gamma} \to P$, given by:
\begin{equation}\label{crnt1}
\begin{split}
& \tau\left(\sigma_P^{\gamma}(f)\right)=\tau\left(c^{\gamma}(f)+N\right)=\tau^{\gamma}(f),\, \text{ for all } f \in P^{\gamma}\, \text{ and } \gamma \in \Omega,\\
& \text{ thus }\tau \circ \sigma_P^{\gamma}=\tau^{\gamma}. 
\end{split}
\end{equation}
Then $\tau$ is a $\Lambda$-isomorphism (see \textbf{Appendix}, \textbf{Prop. A\ref{tauP}}). Observe that $\sigma_A^{\gamma} \circ q^{\gamma}=q \circ \tau^{\gamma}$, for all $\gamma$. This proves that if $f$ belongs to $M^{\gamma}=\Ker(q^{\gamma})$, then $\tau^{\gamma}(f)$ is in $M=\Ker(q)$; that is $\tau^{\gamma}(M^{\gamma}) \subset M$. 
\smallskip

Let $\theta^{\gamma}=\tau^{\gamma}|_{M^{\gamma}}:M^{\gamma} \to M$, then $\theta^{\beta} \circ \chi_{\beta}^{\alpha}=\theta^{\alpha}$ for all $\alpha \leq \beta$. Let $\sigma_M^{\alpha}:M^{\alpha} \to \underrightarrow{\operatorname{lim}}M^{\gamma}$ be the canonical map. There exists a $\Lambda$-homomorphism $\theta:\underrightarrow{\operatorname{lim}}M^{\gamma} \to M$, satisfying:
\begin{equation}\label{M1}
\begin{split}
& \theta\left(\sigma_M^{\gamma}(f)\right)=\theta\left(c^{\gamma}(f)+N\right)=\theta^{\gamma}(f),\, \text{ for all } f \in M^{\gamma}\, \text{ and } \gamma \in \Omega,\\
& \text{ thus } \,\theta \circ \sigma_M^{\gamma}=\theta^{\gamma}.
\end{split}
\end{equation}
Then $\theta$ is an isomorphism (see \textbf{Appendix}, \textbf{Prop. A\ref{thetaP}}). Therefore we have commutative diagram with exact rows:
\begin{equation}\label{d2022-1}
\xymatrix{
0 \ar[r] & M^{\gamma} \ar[d]^{\theta^{\gamma}} \ar[r]^{\iota^{\gamma}} & P^{\gamma} \ar[d]^{\tau^{\gamma}} \ar[r]^{q^{\gamma}} & A^{\gamma} \ar[d]^{\sigma_A^{\gamma}} \ar[r] & 0\\
0 \ar[r] & M \ar[r]^{\iota} & P \ar[r]^{q} & A \ar[r] & 0
}
\end{equation}
Which proves \textbf{Claim 2}. The \textbf{Claim 1} and \textbf{Claim 2} are part of a more general result which affirms that \emph{if $A=\underrightarrow{\operatorname{lim}}A^{\gamma}$, then there exist projective resolutions $ \oplus_k P^{\gamma}_k$ of $A^{\gamma}$, forming a direct system such that $P=\underrightarrow{\operatorname{lim}}P^{\gamma}$ is a projective resolution of $A$} (see \cite{Cartan}, \textbf{Lemma $9.5^{\ast}$, Chapter V, \S 10}). We include the details of the proofs of \textbf{Claim 1} and \textbf{Claim 2} because we need the constructions of the morphisms to obtain our results, as we can see below.
\smallskip

By \eqref{d2022-1}, we have the commutative diagram with exact rows:
$$
\xymatrix{
 T(M^{\gamma}) \ar[d]^{T\left(\theta^{\gamma}\right)} \ar[r]^{T\left(\iota^{\gamma}\right)} & T(P^{\gamma})  \ar[d]^{T\left(\tau^{\gamma}\right)} \ar[r]^{T\left(q^{\gamma}\right)} & T(A^{\gamma}) \ar[d]^{T\left({\sigma_A}^{\gamma}\right)} \ar[r] & 0\\
T(M) \ar[r]^{T(\iota)} & T(P) \ar[r]^{T(q)} & T(A) \ar[r] & 0
}
$$
We also have the direct systems of $\Lambda_1$-modules: $\{T\left(M^{\alpha}\right);T\left(\chi_{\beta}^{\alpha}\right)\}$, $\{T\left(P^{\alpha}\right),T\left(\phi_{\beta}^{\alpha}\right)\}$, and $\{T\left(A^{\alpha}\right);T\left(\varphi_{\beta}^{\alpha}\right)\}$.
\smallskip

Each of the $\Lambda_1$-homomorphisms $T\left(\theta^{\gamma}\right)$, $T\left(\tau^{\gamma}\right)$ and $T\left(\sigma_A^{\gamma}\right)$, for all $\gamma$ in $\Omega$, produce the $\Lambda_1$-homomorphisms:
$$
\begin{array}{rccl}
\widehat{\theta}:& \underrightarrow{\operatorname{lim}}T(M^{\gamma}) & \longrightarrow & T(M)\\
& c^{\gamma}(x)+N & \mapsto & T(\theta^{\gamma})(x)
\end{array}\,\quad 
\begin{array}{rccl}
\widehat{\tau}:& \underrightarrow{\operatorname{lim}}T(P^{\gamma}) & \longrightarrow & T(P)\\
& c^{\gamma}(x)+N & \mapsto & T(\tau^{\gamma})(x)
\end{array}
$$
$$
\begin{array}{rccl}
\sigma^T_A:& \underrightarrow{\operatorname{lim}}T(A^{\gamma}) & \longrightarrow & T(A)\\
& c^{\gamma}(x)+N & \mapsto & T(\sigma_A^{\gamma})(x).
\end{array}
$$
(In each case $N$ is the $\Lambda$-submodule for which the corresponding direct limit is obtained.) We have the commutative diagram with exact rows:
\begin{equation}\label{diagrama-satelite}
\xymatrix{
 \underrightarrow{\operatorname{lim}}T(M^{\gamma}) \ar[d]^{\widehat{\theta}} \ar[r]^{\eta} &  \underrightarrow{\operatorname{lim}}T(P^{\gamma})  \ar[d]^{\widehat{\tau}} \ar[r]^{\zeta} &  \underrightarrow{\operatorname{lim}}T(A^{\gamma}) \ar[d]^{\sigma^T_A}  \ar[r] & 0\\
T(M) \ar[r]^{T(\iota)} & T(P) \ar[r]^{T(q)} & T(A) \ar[r] & 0
}
\end{equation}
where $\eta$ and $\zeta$, are given by:
$$ 
\begin{array}{rccl}
\eta:& \underrightarrow{\operatorname{lim}}T\left(M^{\gamma}\right) & \longrightarrow & \underrightarrow{\operatorname{lim}}T\left(P^{\gamma}\right)\\
& c^{\gamma}(x)+N & \mapsto & c^{\gamma}\left(T(\iota^{\gamma})(x)\right)+N
\end{array}
$$
$$
\begin{array}{rccl}
\zeta:& \underrightarrow{\operatorname{lim}}T\left(P^{\gamma}\right) & \longrightarrow & \underrightarrow{\operatorname{lim}}T\left(A^{\gamma}\right)\\
& c^{\gamma}(x)+N & \mapsto & c^{\gamma}\left(T(q^{\gamma})(x)\right)+N.
\end{array}
$$

\begin{claim}\label{claim 3}{\sl
$\,$

\textbf{(i)} Let $\sigma^T_P:\underrightarrow{\operatorname{lim}}T\left(P^{\gamma}\right) \to T\left(\underrightarrow{\operatorname{lim}}P^{\gamma}\right)$, $c^{\gamma}(x)+N \mapsto T\left(\sigma_P^{\gamma}\right)(x)$ (see \eqref{Lsigma}). Then $T(\tau) \circ \sigma^T_P=\widehat{\tau}$. 
\smallskip

\textbf{(ii)} Let $\sigma^T_M:\underrightarrow{\operatorname{lim}}T\left(M^{\gamma}\right) \to T\left(\underrightarrow{\operatorname{lim}}M^{\gamma}\right)$, $c^{\gamma}(x)+N \mapsto T\left(\sigma_M^{\gamma}\right)(x)$ (see \eqref{Lsigma}). Then $T(\theta) \circ \sigma^T_M=\widehat{\theta}$.} 
\end{claim}
\begin{proof}

\textbf{(i)} By \eqref{crnt1} we have $T(\tau) \circ T\left(\sigma_P^{\gamma}\right)=T\left(\tau^{\gamma}\right)$, for all $\gamma$. Then $T(\tau) \circ \sigma^T_P=\widehat{\tau}$.
\smallskip

\textbf{(ii)} By \eqref{M1}, we get $T(\theta) \circ T\left(\sigma^{\gamma}_M\right)=T\left(\theta^{\gamma}\right)$ for all $\gamma$. Whence $T(\theta) \circ \sigma^T_M=\widehat{\theta}$. Which proves our claim.
\end{proof}

By \textbf{Claim \ref{claim 3}.(i)}, we know that $T(\tau) \circ \sigma^T_P=\widehat{\tau}$. We also know that $\tau$ is an isomorphism (see \textbf{Prop A\ref{tauP}}) and by hypothesis, $\sigma^T_P$ is an isomorphism, whence $\widehat{\tau}$ is an isomorphism. As \eqref{diagrama-satelite} is commutative and it has exact rows, $\sigma^T_A$ is surjective. 
\smallskip

Due to $\{A^{\alpha};\varphi^{\alpha}_{\beta}\}$ is an arbitrary direct system of modules, we can apply the same argument to the direct system of modules $\{M^{\alpha};\chi^{\alpha}_{\beta}\}$, and we get that $\sigma^T_M$ is surjective. Since $\theta$ is an isomorphism (see \textbf{Prop A\ref{thetaP}}) and $\sigma^T_M$ is surjective, from $T(\theta) \circ \sigma^T_M=\widehat{\theta}$ (see \textbf{Claim \ref{claim 3}.(ii)}), it follows that $\widehat{\theta}$ is surjective. 
\smallskip

As \eqref{diagrama-satelite} is commutative and it has exact rows and $\widehat{\tau}$ is an isomorphism, then $\sigma^T_A$ is injective, which proves our result.
\end{proof}

\section*{Appendix}

\begin{prop A}\label{tauP}{\sl
The map $\tau:\underrightarrow{\operatorname{lim}}P^{\gamma} \to P$ is a $\Lambda$-isomorphism.
}
\end{prop A}
\begin{proof}
Let $z$ be in $\underrightarrow{\operatorname{lim}}P^{\gamma}$ such that $\tau(z)=0$. By \textbf{Prop. \ref{rotman}}, there exists $\alpha$ such that $z=c^{\alpha}(f)+N$ where $f$ belongs to $P^{\alpha}$. Let $f=\sum_{a \in A^{\alpha}}f(a)e^{\alpha}(a)$. Then,
\begin{equation}\label{tau1}
0=\tau(z)=\tau^{\alpha}(f)=\sum_{a \in A^{\alpha}}f(a)\tau^{\alpha}\left(e^{\alpha}(a)\right)=\sum_{a \in A^{\alpha}}f(a)e\left(\sigma_A^{\alpha}(a)\right).
\end{equation}
Due to $P$ is a free module generated by $e(x)$, where $x$ is in $A$ it follows from \eqref{tau1} that $f(a)=0$ for all $a$ in $A^{\alpha}$, therefore $f=0$ which proves that $\tau$ is injective. 
\smallskip

Let $e(x)$ be in $P$, where $x$ is in $A=\underrightarrow{\operatorname{lim}}A^{\gamma}$. Then there exists $\alpha$ such that $x=c^{\alpha}(a)+N$, where $a$ is in $A^{\alpha}$. Thus,
\begin{equation}\label{tau2}
\aligned
& \tau\left(c^{\alpha}(e^{\alpha}(a))+N\right)=\tau^{\alpha}\left(e^{\alpha}(a)\right)=e\left( \sigma_A^{\alpha}(a)\right)\\
&=e\left(c^{\alpha}(a)+N\right)=e(x).
\endaligned
\end{equation}
Due to $P$ is a free module generated by $e(x)$, where $x$ is in $A$, then we deduce from \eqref{tau2} that $\tau$ is surjective.  
\end{proof}

\begin{prop A}\label{thetaP}{\sl
The map $\theta:\underrightarrow{\operatorname{lim}}M^{\gamma} \to M$ is a $\Lambda$-isomorphism.}
\end{prop A}
\begin{proof}
As $\theta^{\gamma}$ is the restriction $\tau^{\gamma}_{M^{\gamma}}$, we shall only prove that $\theta$ is surjective. Let $f$ be in $M \subset P$. As $\tau$ is surjective, there exists $h$ in $\underrightarrow{\operatorname{lim}}P^{\gamma}$ such that $f=\tau\left(\sigma_P^{\alpha}(h)\right)$, for some $\alpha$ in $\Omega$. Due to $f$ is in $M=\Ker(q)$, then,
$$
0=q(f)=q(\tau^{\alpha}(h))=\sigma_A^{\alpha}(q^{\alpha}(h))=c^{\alpha}(q^{\alpha}(h))+N.
$$
By \textbf{Prop. \ref{rotman}.(ii)}, there exists $\beta \geq \alpha$ such that $\varphi^{\alpha}_{\beta}\left(q^{\alpha}(h)\right)=q^{\beta}\left(\phi^{\alpha}_{\beta}(h)\right)$, then $\phi^{\alpha}_{\beta}(h)$ belongs to $\Ker(q^{\beta})=M^{\beta}$. In addition observe that:
$$
\tau\left(\sigma^{\beta}_P(\phi^{\alpha}_{\beta}(h))\right)=\tau\left(c^{\beta}(\phi^{\alpha}_{\beta}(h))+N\right)=\tau(c^{\alpha}(h)+N)=f.
$$
Thus $f=\tau^{\beta}\left(\phi^{\alpha}_{\beta}(h)\right)=\theta^{\beta}\left(\phi^{\alpha}_{\beta}(h)\right)=\theta\left(\sigma_M^{\beta}(\phi^{\alpha}_{\beta}(h))\right)$, from it follows our result.
\end{proof}

\section*{Acknowledgement}

The author thanks the support provided by post-doctoral fellowship granted by CONAHCYT 769309.
\smallskip

\textbf{Conflict of interest} The author declares that he has no conflict of interest.

\end{document}